\documentclass[pdflatex,sn-mathphys-ay]{sn-jnl}

\usepackage{graphicx}
\usepackage{amsmath,amssymb,amsfonts}
\usepackage{textcomp}
\usepackage{booktabs}

\AtBeginDocument{\providecommand*{\backrefalt}[4]{}\renewcommand*{\backrefalt}[4]{}}

\theoremstyle{thmstyleone}
\newtheorem{theorem}{Theorem}
\newtheorem{lemma}[theorem]{Lemma}
\newtheorem{proposition}[theorem]{Proposition}
\newtheorem{corollary}[theorem]{Corollary}

\theoremstyle{thmstylethree}
\newtheorem{definition}[theorem]{Definition}
\newtheorem{assumption}[theorem]{Assumption}

\theoremstyle{thmstyletwo}
\newtheorem{remark}{Remark}
\newtheorem*{example}{Example}

\begin{document}

\title[Statistical Properties of Divergence Regularized Optimal Transport]{General Divergence Regularized Optimal Transport: Sample Complexity and Central Limit Theorems}

\author[1]{\fnm{Jiaping} \sur{Yang}}\email{jpyang22@m.fudan.edu.cn}
\author[1]{\fnm{Yunxin} \sur{Zhang}}\email{xyz@fudan.edu.cn}

\affil[1]{\orgdiv{School of Mathematical Sciences}, \orgname{Fudan University},
\orgaddress{\city{Shanghai}, \country{China}}}

\abstract{We study empirical divergence-regularized optimal transport on arbitrary Polish spaces. For bounded continuous costs and dual-regular conjugates $\psi\in C^1(\mathbb R)$, we establish a dimension-free $n^{-1/2}$ bound for the empirical transport value using dual interpolation, Hoeffding projection, and bounded-difference arguments. For the asymptotic theory, we isolate the population identifiability condition needed for stability and give sufficient conditions for population uniqueness. In particular, active-set identifiability implies uniqueness of the canonical population potentials, and connectedness of either population support is sufficient. This condition requires neither compactness nor Euclidean structure and covers quadratic regularization, whose conjugate is $C^1$ but not $C^2$. We also give finite-support and large-regularization criteria. Under a bounded Lipschitz cost, the empirical potentials are stable, the transport value is asymptotically linear, and one- and two-sample central limit theorems hold with population centering. This includes noncompact and infinite-dimensional Polish state spaces.}


\maketitle

\section{Introduction}
Optimal transport (OT) is widely used to compare probability distributions, including in statistical problems involving empirical distributions. When the underlying marginals are unknown and replaced by empirical measures, the resulting transport functionals raise two basic questions: how accurately they estimate their population counterparts at finite sample sizes, and whether they admit tractable asymptotic distributions. These questions arise, for example, when OT is used as a statistical discrepancy or as an ingredient in statistical inference; see \cite{MR3939527,MR4659653,peyre2019computational}.

Regularization is often introduced to improve the computational and statistical behavior of empirical OT. Entropic regularization is especially prominent because it leads to efficient Sinkhorn-type algorithms with well-developed convergence theory \citep{cuturi2013sinkhorn,MR4418040,MR4506579,MR4546630}. Its optimal couplings are, however, dense, and the constants governing computation and statistical approximation may deteriorate as the regularization parameter decreases. Divergence-based alternatives provide a flexible way to retain regularization while allowing sparse optimal couplings.

This paper studies the empirical value of divergence-regularized OT. We focus on dimension-free finite-sample accuracy and Gaussian limits. The two problems depend on different aspects of the dual formulation. Finite-sample control can be obtained from bounded dual geometry, whereas a first-order limit theory requires stability of the population dual representation. For sparsity-inducing regularizers this distinction is nontrivial, because the dual potentials need not be unique. Accordingly, the limit theory is preceded by an analysis of uniqueness of the population dual potentials.

\subsection{Related work}
Non-entropic regularization was introduced in \cite{blondel2018smooth} through quadratic regularization and extended to broader classes of regularizers in \cite{MR3862422}. The resulting divergence-regularized problems are supported by a general duality theory \citep{MR4409806}, structural analyses of sparse optimal plans and vanishing regularization \citep{wiesel2024,garriz2024infinitesimal}, approximation bounds \citep{MR4755769}, and numerical methods \citep{lorenz2019orlicz,pasechnyuk}.

Classical empirical OT can exhibit dimension-dependent convergence rates; see \cite{MR4901218} for an overview. For entropic regularization, parametric bounds for the transport value were obtained under compact support in \cite{genevay2019sample} and under broader tail conditions in \cite{mena2019statistical}; weak limits and related inferential questions have been studied in \cite{gonzalez2022weak,gonzalez2023weak,pooladian2021entropic,MR4105566,MR4616887}. For general divergences, \citet{MR4815982} established stability and sample-complexity bounds involving smoothness and intrinsic dimension. More recently, \citet{GonzalezSanzEcksteinNutz2025} developed a parametric theory for sparse regularized transport on compact subsets of $\mathbb R^d$. Under a connected-support condition, a $C^2$ conjugate, and a $C^1$ cost for the central limit theorems, they obtained Gaussian limits for the transport value, the dual potentials, and the optimal coupling.

Here we restrict attention to the scalar transport value but work on arbitrary Polish state spaces, without compact-support or finite-dimensional Euclidean assumptions. The finite-sample theory requires a bounded continuous cost and a dual-regular $C^1$ conjugate, whereas the limit theory assumes a bounded Lipschitz cost under the same dual regularity. This setting includes quadratic regularization, whose conjugate is $C^1$ but not $C^2$. Population uniqueness is characterized through active-set identifiability, with connectedness of either support as a sufficient condition. Thus our value results complement the potential and coupling limit theory of \cite{GonzalezSanzEcksteinNutz2025}.

\subsection{Contributions and organization}
The paper has three main results. First, we establish a dimension-free $n^{-1/2}$ bound for the empirical divergence-regularized value under bounded continuous costs and dual regularity on arbitrary Polish spaces. The proof combines dual interpolation, projection onto additive functions, and bounded differences, without differentiating the cost or imposing compact support.

Second, we give sufficient conditions for population identifiability of the dual potentials. Active-set identifiability implies uniqueness, and connectedness of either population support is sufficient on arbitrary Polish spaces. For finitely supported marginals, the condition reduces to connectedness of an active bipartite graph. We also obtain a sufficient condition in the large-$\epsilon$ regime.

Third, under bounded Lipschitz costs and population identifiability, we prove local-uniform stability of canonical empirical potentials and derive an asymptotically linear representation of the transport value. This yields one- and two-sample central limit theorems. The same stability argument makes the bias negligible on the central-limit scale, so the limiting distributions may be centered at the population value. Together with the connected-support criterion, the result applies to noncompact and infinite-dimensional Polish spaces.

Section \ref{Se2} introduces the notation and the duality and regularity facts used throughout the paper. Section \ref{Se3} establishes the finite-sample bounds. Section \ref{Se4} develops the identifiability criteria and the central limit theory. Section \ref{Se5} concludes, and the proof of the projection lemma is deferred to the \hyperref[app]{Appendix}.

\section{Preliminaries}\label{Se2}

\subsection{Notation}\label{13}
Throughout, $(\mathcal X,d_x)$ and $(\mathcal Y,d_y)$ are complete separable metric spaces, and $\mu$ and $\nu$ are Borel probability measures on them. For measures $\mu$ and $\eta$ on the same measurable space, $\mu\ll\eta$ denotes absolute continuity and $\mathrm d\mu/\mathrm d\eta$ the corresponding Radon--Nikodym derivative. We write $\operatorname{supp}(\mu)$ for the support of $\mu$, $\mu\otimes\nu$ for the product measure, and $\Pi(\mu,\nu)$ for the set of couplings of $\mu$ and $\nu$. Given independent samples $X_1,\ldots,X_n\sim\mu$ and $Y_1,\ldots,Y_m\sim\nu$, let $\mu_n$ and $\nu_m$ denote the associated empirical measures. We use $\mathrm E_\mu$, $\mathrm{Var}_\mu$, and $\mathrm{Cov}_\mu$ for expectation, variance, and covariance under $\mu$, omitting the subscript when the underlying measure is clear. Finally, $\mathcal N(0,1)$ denotes the standard normal law and $\stackrel{D}{\longrightarrow}$ convergence in distribution.

For functions $f:\mathcal X\to\mathbb R$ and $g:\mathcal Y\to\mathbb R$, we write $\|f\|_\infty=\sup_x|f(x)|$ and $\operatorname{osc}(f)=\sup f-\inf f$. The notation $h=f\oplus g$ means $h(x,y)=f(x)+g(y)$, and $\mathcal B_b^\oplus$ denotes bounded measurable additive functions of this form. For a probability measure $P$, $Pf$ abbreviates $\int f\,\mathrm dP$; similarly, for $k(x,y)$ we write $Qk(x)=\int k(x,y)Q(\mathrm dy)$ and $Pk(y)=\int k(x,y)P(\mathrm dx)$. Inner products and norms are always taken in the indicated $L^2$ space. We abbreviate $x_+=\max(x,0)$ and use $C$ or $K$ for finite positive constants that may change from line to line. For nonnegative sequences, $a_n\lesssim b_n$ means $a_n\le Cb_n$ with $C$ independent of $n$.

\subsection{Dual formulation and regularity}
We use the following dual formulation and regularity assumptions.

Let $\phi:[0,\infty)\to\mathbb R\cup\{+\infty\}$ be proper, lower semicontinuous, strictly convex, and superlinear, with $\phi(1)=0$. Its one-sided convex conjugate is
\begin{equation*}
 \psi(t)=\sup_{r\ge 0}\{rt-\phi(r)\},\qquad t\in\mathbb R.
\end{equation*}
For a coupling $\pi\in\Pi(\mu,\nu)$, define
\begin{equation*}
 \mathrm H_\phi(\pi\|\mu\otimes\nu)=
 \begin{cases}
 \displaystyle\int \phi\!\left(\frac{\mathrm d\pi}{\mathrm d(\mu\otimes\nu)}\right)\mathrm d(\mu\otimes\nu),&\pi\ll\mu\otimes\nu,\\[1ex]
 +\infty,&\text{otherwise}.
 \end{cases}
\end{equation*}
For $\epsilon>0$, the divergence-regularized transport value is
\begin{equation}\label{E1}
 S_\epsilon(\mu,\nu)=\inf_{\pi\in\Pi(\mu,\nu)}\left\{\int c\,\mathrm d\pi+\epsilon\,\mathrm H_\phi(\pi\|\mu\otimes\nu)\right\}.
\end{equation}
The product coupling is feasible and has finite divergence, so the problem is proper. Under the standing assumptions below, the primal objective is lower semicontinuous on the weakly compact set $\Pi(\mu,\nu)$; hence the infimum is attained. Strict convexity of $\phi$ yields uniqueness of the optimal coupling, although the dual potentials need not be unique.

\begin{definition}[Dual regularity]\label{def:dualregular}
The divergence generator $\phi$ is called \emph{dual regular} if $\psi\in C^1(\mathbb R)$ and there exist $t_0\in\mathbb R$ and $\delta,C_0>0$ such that
\begin{enumerate}
 \item $\psi$ is strictly convex on $[t_0-\delta,\infty)$ and $\psi'(t_0)=1$;
 \item $\psi'(t)\ge t$ for all $t\ge C_0$.
\end{enumerate}
\end{definition}
This is the dual-regularity condition used in \cite{MR4815982}. Notice in particular that the normalization is $\psi'(t_0)=1$, rather than $\psi(t_0)=1$.

\begin{assumption}[Standing assumptions]\label{ass:base}
The cost satisfies $c\in C_b(\mathcal X\times\mathcal Y)$, and $\phi$ is dual regular in the sense of Definition \ref{def:dualregular}.
\end{assumption}

For $\epsilon=1$, the dual functional can be written entirely in terms of the additive sum $h=f\oplus g$, where $h(x,y)=f(x)+g(y)$:
\begin{equation}\label{EE1}
 L_{\alpha,\beta}(h)
 :=\int\big\{h(x,y)-\psi(h(x,y)-c(x,y))\big\}\,\mathrm d(\alpha\otimes\beta)(x,y).
\end{equation}
Here $\alpha,\beta$ denote arbitrary probability measures on $\mathcal X,\mathcal Y$. We write
\begin{equation*}
 \ell_h(x,y):=h(x,y)-\psi(h(x,y)-c(x,y)).
\end{equation*}
The additive representation is gauge invariant: replacing $(f,g)$ by $(f+a,g-a)$ leaves both $h$ and the dual objective unchanged.

We use the following consequences of the dual theory in \cite{MR4815982}.

\begin{proposition}[Dual attainment and regularity]\label{prop:dualbasics}
Under Assumption \ref{ass:base}, for every pair of probability measures $(\alpha,\beta)$ on $(\mathcal X,\mathcal Y)$ the following hold.
\begin{enumerate}
 \item The primal optimizer is unique, strong duality holds, and the dual supremum is attained:
 \begin{equation}\label{eq:dualbounded}
 S(\alpha,\beta)=\sup_{h\in\mathcal B_b^\oplus}L_{\alpha,\beta}(h).
 \end{equation}
 \item Every dual-optimal equivalence class admits a canonical representative on all of $\mathcal X\times\mathcal Y$ satisfying the marginal equations
 \begin{align}
 \int \psi'(f(x)+g(y)-c(x,y))\,\beta(\mathrm dy)&=1, &&x\in\mathcal X,\label{eq:margx}\\
 \int \psi'(f(x)+g(y)-c(x,y))\,\alpha(\mathrm dx)&=1, &&y\in\mathcal Y.\label{eq:margy}
 \end{align}
 The corresponding optimal coupling is
 \begin{equation}\label{eq:recovery}
 \frac{\mathrm d\pi^*}{\mathrm d(\alpha\otimes\beta)}(x,y)=\psi'(f(x)+g(y)-c(x,y)).
 \end{equation}
 \item There is a deterministic constant $B<\infty$, depending only on $\|c\|_\infty$ and $\phi$, such that every canonical optimizer, after the normalization $\int f\,\mathrm d\alpha=0$, satisfies
 \begin{equation}\label{eq:uniformpotbound}
 \|f\|_\infty+\|g\|_\infty+\|f\oplus g\|_\infty\le B.
 \end{equation}
 Accordingly, with
 \begin{equation*}
 P_B:=\sup_{|t|\le B+\|c\|_\infty}|\psi'(t)|,
 \qquad
 M_B:=B+\sup_{|t|\le B+\|c\|_\infty}|\psi(t)|,
 \end{equation*}
 both the optimal density in \eqref{eq:recovery} and the dual integrand $\ell_h$ are bounded by deterministic constants independent of $(\alpha,\beta)$.
 \item Every canonical extension is continuous. More precisely,
 \begin{align}
  |f(x)-f(x')|&\le \sup_{y\in\operatorname{supp}(\beta)}|c(x,y)-c(x',y)|,\label{eq:directionalf}\\
  |g(y)-g(y')|&\le \sup_{x\in\operatorname{supp}(\alpha)}|c(x,y)-c(x,y')|.\label{eq:directionalg}
 \end{align}
 In particular, if $c$ is $L$-Lipschitz in each variable, then the canonical extensions of $f$ and $g$ are $L$-Lipschitz.
\end{enumerate}
\end{proposition}
\begin{proof}
Lemma 2.1 of \cite{MR4815982} gives primal attainment and strong duality, and Lemma 2.3 there gives dual attainment together with the pointwise marginal equations. Strict convexity of $\phi$ on densities and convexity of the coupling set imply uniqueness of the primal optimizer. The canonical representatives and uniform bounds follow from the bounded-cost construction in the proof of the latter lemma. Two details will be used repeatedly. First, the scalar extension equations are normalized by $\psi'(t_0)=1$ from Definition \ref{def:dualregular}, which determines the canonical values off the supports. Second, the oscillation bounds hold for every canonical optimizer and are controlled by the oscillation of $c$. After fixing the gauge by $\int f\,\mathrm d\alpha=0$, these estimates give a deterministic bound for $f$. For each $y$, the second marginal equation and $\psi'(t_0)=1$ place $g(y)$ between
\begin{equation*}
 t_0-\sup_x\{f(x)-c(x,y)\}
 \quad\text{and}\quad
 t_0-\inf_x\{f(x)-c(x,y)\},
\end{equation*}
so $g$ is bounded by the same standing quantities. This yields a bound of the form \eqref{eq:uniformpotbound} depending only on $\|c\|_\infty$ and $\phi$. Since $h-c$ therefore ranges over a fixed compact interval, continuity of $\psi$ and $\psi'$ gives the final boundedness statement.

The continuity and directional bounds follow directly from the scalar marginal equations. Fix a bounded canonical $g$ and define
\begin{equation*}
 F_x(s):=\int\psi'(s+g(y)-c(x,y))\,\beta(\mathrm dy).
\end{equation*}
Dual regularity makes the level-one solution of $F_x(s)=1$ unique. The map $F_x$ is nondecreasing. At any solution, the set on which $\psi'(s+g-c)>\psi'(t_0-\delta)$ must have positive $\beta$-mass; otherwise $F_x(s)\le\psi'(t_0-\delta)<1$. On this set the argument lies above $t_0-\delta$, where $\psi'$ is strictly increasing. Hence every positive increment of $s$ increases $F_x(s)$ strictly, excluding two distinct level-one solutions.

If $x_n\to x$, uniform boundedness of the potentials and boundedness of $c$ confine all arguments of $\psi'$ to a fixed compact interval. Since $c(x_n,y)\to c(x,y)$ pointwise, dominated convergence yields $F_{x_n}(s)\to F_x(s)$ for each fixed $s$. Along any subsequence for which $f(x_{n_k})\to s_*$, the same argument gives $F_x(s_*)=\lim_kF_{x_{n_k}}(f(x_{n_k}))=1$. Uniqueness of the level-one solution then gives $s_*=f(x)$. Hence $f(x_n)\to f(x)$, so $f$ is continuous; the proof for $g$ is symmetric.

For the directional estimate, set
\begin{equation*}
 D(x,x'):=\sup_{y\in\operatorname{supp}(\beta)}|c(x,y)-c(x',y)|.
\end{equation*}
Then $c(x,y)\le c(x',y)+D(x,x')$ on $\operatorname{supp}(\beta)$, so monotonicity of $\psi'$ and \eqref{eq:margx} give
\begin{equation*}
 \int\psi'\big(f(x')+D(x,x')+g(y)-c(x,y)\big)\,\beta(\mathrm dy)\ge1.
\end{equation*}
Uniqueness of the level-one scalar solution yields $f(x)\le f(x')+D(x,x')$. Interchanging $x$ and $x'$ proves \eqref{eq:directionalf}; \eqref{eq:directionalg} is symmetric. The Lipschitz statement follows immediately. No uniqueness of the optimal pair $(f,g)$ is used here.
\end{proof}

\begin{remark}\label{rem:epsscale}
It suffices to prove the statistical results for $\epsilon=1$. Indeed, if $\bar c=c/\epsilon$ and $S_1^{\bar c}$ denotes the problem with unit regularization parameter and cost $\bar c$, then
\begin{equation*}
 S_\epsilon^c(\mu,\nu)=\epsilon S_1^{\bar c}(\mu,\nu).
\end{equation*}
If $(\bar f,\bar g)$ is optimal for the unit-regularization problem with cost $c/\epsilon$, then $(\epsilon\bar f,\epsilon\bar g)$ is optimal for $S_\epsilon^c$. Hence the rates below remain unchanged for fixed $\epsilon>0$, although their constants depend on $\epsilon$ through the rescaled cost $c/\epsilon$ and the associated potential bounds. The dependence of these constants on $\epsilon$ is not optimized here.
\end{remark}

In the remainder of the paper we set $\epsilon=1$ and suppress it from the notation.

Before studying sampling error, we record the measurability needed for the expectations and concentration inequalities below.

\begin{lemma}\label{lem:valuemeas}
Under Assumption \ref{ass:base}, for every $n,m\ge1$ the maps
\begin{equation*}
 (X_1,\ldots,X_n,Y_1,\ldots,Y_m)\longmapsto S(\mu_n,\nu_m),
\end{equation*}
and the corresponding one-sample maps defining $S(\mu_n,\nu)$ and $S(\mu,\nu_m)$ are Borel measurable.
\end{lemma}
\begin{proof}
For the two-sample problem, partition the sample space according to the finitely many equality patterns among the $X_i$'s and among the $Y_j$'s. On each resulting Borel stratum, Proposition \ref{prop:dualbasics} reduces the dual value to the maximum of a function that is continuous in finitely many potential values and Borel measurable in the sample, over a fixed compact set after gauge normalization and the uniform bound \eqref{eq:uniformpotbound}. The measurable maximum theorem therefore makes the value Borel measurable on each stratum, and hence on the whole sample space. For $S(\mu_n,\nu)$, one first eliminates the $g$-coordinate pointwise through the unique scalar marginal equation and then maximizes the resulting finite-dimensional reduced objective in the empirical support values of $f$; this reduced objective is continuous in those values and Borel measurable in the sample. The same measurable-maximum argument applies. The case $S(\mu,\nu_m)$ is symmetric. A detailed version of this construction, including measurable canonical optimizers, is given in Lemma \ref{lem:measselection} below.
\end{proof}

\section{Finite-sample bounds for the regularized transport value}\label{Se3}
The finite-sample analysis begins with a dimension-free bound for the empirical regularized transport value. The argument differs from the empirical-process approach in \cite{MR4815982}: dual interpolation reduces the bias to first-order Hoeffding projections, while bounded differences control the variance. No derivatives of the cost or higher-order smoothness of the potentials are required.

Let
\begin{equation*}
 S_{n,m}:=S(\mu_n,\nu_m),\qquad S:=S(\mu,\nu),
\end{equation*}
where the two samples are independent. The principal equal-sample statement is the following.

\begin{theorem}[Dimension-free sample complexity]\label{T5}
Under Assumption \ref{ass:base}, there exists a constant $C<\infty$, depending only on the bounded-cost and regularization constants, such that
\begin{equation*}
 \mathrm E\big|S(\mu_n,\nu_n)-S(\mu,\nu)\big|\le \frac{C}{\sqrt n}.
\end{equation*}
\end{theorem}

The equal-sample theorem follows from more general $n,m$ bounds for the bias and stochastic fluctuation. Two auxiliary lemmas isolate the mechanism behind the bias estimate.

\begin{lemma}[Dual interpolation]\label{lem:interpolation}
Let $\alpha,\beta$ be probability measures and let $h_0,h_1\in\mathcal B_b^\oplus$. Define $h_t=(1-t)h_0+th_1$ and $\varphi(t)=L_{\alpha,\beta}(h_t)$. Then $\varphi$ is concave and continuously differentiable, with
\begin{equation}\label{ET1}
 \varphi'(t)=\int (h_1-h_0)\big[1-\psi'(h_t-c)\big]\,\mathrm d(\alpha\otimes\beta).
\end{equation}
If $h_0$ is optimal for $(\alpha,\beta)$, then
\begin{equation}\label{eq:interpbound}
 0\le L_{\alpha,\beta}(h_0)-L_{\alpha,\beta}(h_1)
 \le \int (h_1-h_0)\big[\psi'(h_1-c)-1\big]\,\mathrm d(\alpha\otimes\beta).
\end{equation}
\end{lemma}
\begin{proof}
The map $t\mapsto h_t$ is affine and $\psi$ is convex, so $t\mapsto-\psi(h_t-c)$ is concave. Differentiability follows from $\psi\in C^1$ and boundedness of the interpolation segment, giving \eqref{ET1}. Since $\varphi'$ is nonincreasing,
\begin{equation*}
 \varphi(0)-\varphi(1)=-\int_0^1\varphi'(t)\,\mathrm dt\le-\varphi'(1),
\end{equation*}
which is precisely \eqref{eq:interpbound}.
\end{proof}

\begin{lemma}[Additive projection]\label{L31}
Let $P,Q$ be probability measures, $k\in L^2(P\otimes Q)$, and $h=f\oplus g\in L^2(P\otimes Q)$. Then
\begin{equation}\label{eq:additiveprojection}
 |\langle h,k\rangle_{L^2(P\otimes Q)}|^2
 \le \|h\|_{L^2(P\otimes Q)}^2
 \left(
 \|Qk\|_{L^2(P)}^2+\|Pk\|_{L^2(Q)}^2
 \right),
\end{equation}
where $Qk(x)=\int k(x,y)Q(\mathrm dy)$ and $Pk(y)=\int k(x,y)P(\mathrm dx)$.
\end{lemma}
The Appendix gives the elementary Hoeffding-projection proof. The projection is useful here because $p^*-1$ has vanishing population marginal projections; empirical sampling therefore leaves only the corresponding first-order projection errors.

\begin{theorem}[Bias bound]\label{KK1}
Under Assumption \ref{ass:base}, for all $n,m\ge1$,
\begin{equation}\label{eq:biasnm}
 0\le \mathrm E S(\mu_n,\nu_m)-S(\mu,\nu)
 \le C\sqrt{\frac1n+\frac1m}.
\end{equation}
The one-sample analogues also hold:
\begin{equation}\label{eq:biasonesample}
 0\le \mathrm E S(\mu_n,\nu)-S(\mu,\nu)\le \frac C{\sqrt n},
 \qquad
 0\le \mathrm E S(\mu,\nu_m)-S(\mu,\nu)\le \frac C{\sqrt m}.
\end{equation}
\end{theorem}
\begin{proof}
Let $h^*$ be a canonical population optimizer and $h_{n,m}$ any canonical optimizer for $(\mu_n,\nu_m)$. Set
\begin{equation*}
 p^*(x,y):=\psi'(h^*(x,y)-c(x,y)).
\end{equation*}
Proposition \ref{prop:dualbasics} gives uniform bounds for $p^*$ and both optimal sums, while the marginal equations imply
\begin{equation}\label{eq:pcentering}
 \nu(p^*-1)(x)=0\quad\text{for all }x,
 \qquad
 \mu(p^*-1)(y)=0\quad\text{for all }y.
\end{equation}
Since $h_{n,m}$ maximizes $L_{\mu_n,\nu_m}$, Lemma \ref{lem:interpolation} gives
\begin{equation}\label{eq:biaspointwise}
 0\le S_{n,m}-L_{\mu_n,\nu_m}(h^*)
 \le \langle h^*-h_{n,m},p^*-1\rangle_{L^2(\mu_n\otimes\nu_m)}.
\end{equation}
Applying Lemma \ref{L31} and the deterministic bound on $h^*-h_{n,m}$ yields
\begin{equation*}
 S_{n,m}-L_{\mu_n,\nu_m}(h^*)
 \le C\left(
 \|\nu_m(p^*-1)\|_{L^2(\mu_n)}^2
 +\|\mu_n(p^*-1)\|_{L^2(\nu_m)}^2
 \right)^{1/2}.
\end{equation*}
Using \eqref{eq:pcentering}, independence of the samples, and boundedness of $p^*$,
\begin{align*}
 \mathrm E\|\nu_m(p^*-1)\|_{L^2(\mu_n)}^2&\le \frac{C}{m},\\
 \mathrm E\|\mu_n(p^*-1)\|_{L^2(\nu_m)}^2&\le \frac{C}{n}.
\end{align*}
Taking expectations and applying Cauchy--Schwarz proves the upper bound in \eqref{eq:biasnm}. The lower bound follows from
\begin{equation*}
 \mathrm E S_{n,m}=\mathrm E\sup_hL_{\mu_n,\nu_m}(h)
 \ge \mathrm E L_{\mu_n,\nu_m}(h^*)=L_{\mu,\nu}(h^*)=S.
\end{equation*}
For the one-sample problem $(\mu_n,\nu)$, the term $\nu(p^*-1)$ vanishes identically and only the second empirical projection remains; this proves the first inequality in \eqref{eq:biasonesample}. The second one-sample bound follows symmetrically.
\end{proof}

The remaining stochastic fluctuation is controlled by a bounded-difference argument.

\begin{theorem}[Variance bound]\label{V1}
Under Assumption \ref{ass:base},
\begin{equation}\label{eq:varnm}
 \operatorname{Var}(S(\mu_n,\nu_m))\le C\left(\frac1n+\frac1m\right).
\end{equation}
In particular, $\operatorname{Var}(S(\mu_n,\nu_n))\le C/n$.
\end{theorem}
\begin{proof}
Let $h_{n,m}$ be a canonical optimizer and recall $\ell_h=h-\psi(h-c)$. Proposition \ref{prop:dualbasics} gives a deterministic $M<\infty$ such that $|\ell_{h_{n,m}}|\le M$ for all samples. Replace $X_1$ by an independent copy $X_1'$ and denote the resulting empirical measure by $\mu_n'$. Optimality of $h_{n,m}$ gives
\begin{align*}
 S(\mu_n,\nu_m)-S(\mu_n',\nu_m)
 &\le (\mu_n-\mu_n')\nu_m\ell_{h_{n,m}}\\
 &\le \frac{2M}{n}.
\end{align*}
Interchanging the original and primed samples and using a canonical optimizer for the primed problem gives the reverse inequality. Thus replacing any $X_i$ changes the value by at most $2M/n$, and replacing any $Y_j$ changes it by at most $2M/m$. Applying the Efron--Stein inequality in the form $\operatorname{Var}(Z)\le \frac12\sum_i\mathrm E(Z-Z^{(i)})^2$ \citep{MR3185193} gives
\begin{equation*}
 \operatorname{Var}(S_{n,m})
 \le \frac12\left[n\left(\frac{2M}{n}\right)^2+m\left(\frac{2M}{m}\right)^2\right],
\end{equation*}
which proves the $O(n^{-1}+m^{-1})$ variance bound in \eqref{eq:varnm}.
\end{proof}

\begin{proof}[Proof of Theorem \ref{T5}]
Cauchy--Schwarz and the bias--variance decomposition give
\begin{align*}
 \mathrm E|S(\mu_n,\nu_n)-S|
 &\le \left(\mathrm E|S(\mu_n,\nu_n)-S|^2\right)^{1/2}\\
 &=\left(\operatorname{Var}(S(\mu_n,\nu_n))
 +|\mathrm ES(\mu_n,\nu_n)-S|^2\right)^{1/2}.
\end{align*}
Theorems \ref{KK1} and \ref{V1} give the claimed $n^{-1/2}$ rate.
\end{proof}

The same sensitivity estimate also yields concentration around the mean.

\begin{corollary}[Concentration bound]\label{cor:highprob}
Under Assumption \ref{ass:base}, there are constants $C_0,C_1<\infty$ such that for every $t>0$, with probability at least $1-2e^{-t}$,
\begin{equation}\label{eq:highprob}
 |S(\mu_n,\nu_n)-S(\mu,\nu)|
 \le \frac{C_0+C_1\sqrt t}{\sqrt n}.
\end{equation}
Consequently, for $t\ge1$ the right-hand side is bounded by $C\sqrt{t/n}$.
\end{corollary}
\begin{proof}
The sensitivity bound in Theorem \ref{V1} is $2M/n$ for each of the $2n$ observations. McDiarmid's inequality therefore yields
\begin{equation*}
 \mathbb P\left(|S(\mu_n,\nu_n)-\mathrm ES(\mu_n,\nu_n)|>2M\sqrt{t/n}\right)\le2e^{-t}.
\end{equation*}
Combining this estimate with Theorem \ref{KK1} proves the claim.
\end{proof}

\begin{remark}
For a general fixed regularization parameter $\epsilon>0$, Remark \ref{rem:epsscale} applies with the rescaled cost $c/\epsilon$. The rate remains parametric, while the constant is determined by the bounded dual interval generated by $\|c\|_\infty/\epsilon$ and $\psi$. The dependence of this constant on $\epsilon$ is not optimized here.
\end{remark}

\begin{remark}
The variance is of order $n^{-1}$ by bounded differences, whereas the $n^{-1/2}$ bias is governed by the first-order empirical projections of the population optimal density. The bound \eqref{eq:biaspointwise} improves to little-$o$ on the central-limit scale once $\|h_{n,m}-h^*\|$ converges to zero; this stability is established in Section \ref{Se4} under population identifiability.
\end{remark}

\section{Identifiability and central limit theory}\label{Se4}
Distributional limits require one additional ingredient: identifiability of the \emph{population} dual potentials. For sparsity-inducing divergences, differentiability---even $C^2$ regularity of $\psi$---does not by itself imply uniqueness. We therefore isolate the population uniqueness condition used by the stability argument and then give structural sufficient conditions.

\begin{assumption}[Population uniqueness]\label{ass:unique}
The canonical population optimal potential pair is unique up to the additive gauge. Equivalently, after the normalization
\begin{equation}\label{eq:populationnormalization}
 \int f^*\,\mathrm d\mu=0,
\end{equation}
there is a unique canonical pair $(f^*,g^*)$ solving the population dual problem and the marginal equations \eqref{eq:margx}--\eqref{eq:margy}.
\end{assumption}

Assumption \ref{ass:unique} can be checked through the region in which the dual regularizer is strictly curved. Set
\begin{equation}\label{eq:rdelta}
 r_\delta:=\psi'(t_0-\delta)<\psi'(t_0)=1,
\end{equation}
write $S_X:=\operatorname{supp}(\mu)$ and $S_Y:=\operatorname{supp}(\nu)$, and choose any canonical population optimizer $(f^*,g^*)$. By uniqueness of the primal optimizer, the function
\begin{equation}\label{eq:rhostar}
 \rho^*(x,y):=\psi'(f^*(x)+g^*(y)-c(x,y)),\qquad (x,y)\in S_X\times S_Y,
\end{equation}
is independent of the chosen dual optimizer. Define the \emph{strong active set}
\begin{equation}\label{eq:activeset}
 \mathcal A^*:=\{(x,y)\in S_X\times S_Y:\rho^*(x,y)>r_\delta\}.
\end{equation}
On $\mathcal A^*$, every dual argument lies in $(t_0-\delta,\infty)$, where $\psi'$ is strictly increasing.

\begin{definition}[Active-set identifiability]\label{def:activeid}
The population problem is said to satisfy \emph{active-set identifiability} if the following implication holds for bounded continuous functions $u:S_X\to\mathbb R$ and $v:S_Y\to\mathbb R$:
\begin{equation}\label{eq:activeid}
 u(x)+v(y)=0\quad\text{for all }(x,y)\in\mathcal A^*
 \quad\Longrightarrow\quad
 u\equiv a,\quad v\equiv-a
\end{equation}
for some $a\in\mathbb R$.
\end{definition}

\begin{proposition}[Active-set uniqueness]\label{prop:activeunique}
Under Assumption \ref{ass:base}, active-set identifiability implies Assumption \ref{ass:unique}.
\end{proposition}
\begin{proof}
Let $(f_1,g_1)$ and $(f_2,g_2)$ be canonical population optimal pairs. Strict convexity of $\phi$ makes the primal optimizer unique; hence Proposition \ref{prop:dualbasics} gives
\begin{equation*}
 \psi'(f_1\oplus g_1-c)=\psi'(f_2\oplus g_2-c)
 \qquad \mu\otimes\nu\text{-a.e.}
\end{equation*}
Both sides are continuous by Proposition \ref{prop:dualbasics}. Because $\operatorname{supp}(\mu\otimes\nu)=S_X\times S_Y$, the almost-sure equality extends to every point of $S_X\times S_Y$. On $\mathcal A^*$ the common value exceeds $r_\delta$, so monotonicity places both arguments above $t_0-\delta$; strict increase of $\psi'$ on $[t_0-\delta,\infty)$ then gives
\begin{equation*}
 (f_1-f_2)(x)+(g_1-g_2)(y)=0,
 \qquad (x,y)\in\mathcal A^*.
\end{equation*}
Definition \ref{def:activeid} therefore yields $f_1-f_2\equiv a$ on $S_X$ and $g_1-g_2\equiv-a$ on $S_Y$. Finally, the pointwise scalar equations \eqref{eq:margx}--\eqref{eq:margy} uniquely determine the canonical values off the supports, so the same gauge relation holds on all of $\mathcal X$ and $\mathcal Y$. The normalization \eqref{eq:populationnormalization} fixes $a=0$.
\end{proof}

\begin{theorem}[Connected-support uniqueness]\label{thm:connectedunique}
Under Assumption \ref{ass:base}, if either $S_X$ or $S_Y$ is connected, then active-set identifiability holds and hence Assumption \ref{ass:unique} is automatic.
\end{theorem}
\begin{proof}
Assume that $S_X$ is connected; the case of connected $S_Y$ is symmetric. Let bounded continuous $u,v$ satisfy $u\oplus v=0$ on $\mathcal A^*$. For each $x\in S_X$, the first marginal equation gives
\begin{equation*}
 \int\rho^*(x,y)\,\nu(\mathrm dy)=1.
\end{equation*}
Because $r_\delta<1$, the section $\{y\in S_Y:(x,y)\in\mathcal A^*\}$ is nonempty. Choose $y_x$ in this section. Continuity of $\rho^*$ makes $\mathcal A^*$ relatively open in $S_X\times S_Y$; hence some relative neighborhood $U_x$ of $x$ satisfies $(z,y_x)\in\mathcal A^*$ for all $z\in U_x$. Therefore
\begin{equation*}
 u(z)=-v(y_x)=u(x),\qquad z\in U_x.
\end{equation*}
Thus $u$ is locally constant on the connected space $S_X$, and hence $u\equiv a$. For each $y\in S_Y$, the second marginal equation provides some $x_y\in S_X$ with $(x_y,y)\in\mathcal A^*$, so $v(y)=-u(x_y)=-a$. This proves active-set identifiability, and Proposition \ref{prop:activeunique} yields uniqueness.
\end{proof}

\begin{corollary}[Active-graph criterion]\label{cor:activegraph}
Suppose
\begin{equation*}
 \mu=\sum_{i=1}^I p_i\delta_{x_i},\qquad
 \nu=\sum_{j=1}^J q_j\delta_{y_j},
\end{equation*}
where $x_1,\ldots,x_I$ and $y_1,\ldots,y_J$ are the distinct support points and all weights are positive. Form the bipartite graph with vertices $\{x_i\}_{i=1}^I\sqcup\{y_j\}_{j=1}^J$ and an edge $(x_i,y_j)$ whenever $\rho^*(x_i,y_j)>r_\delta$. If this active graph is connected, then Assumption \ref{ass:unique} holds.
\end{corollary}
\begin{proof}
If $u(x_i)+v(y_j)=0$ on every active edge, propagation along paths in the connected bipartite graph shows that all $u(x_i)$ equal some $a$ and all $v(y_j)$ equal $-a$. Definition \ref{def:activeid} therefore holds, and Proposition \ref{prop:activeunique} applies.
\end{proof}

Connectivity is one route to identifiability; curvature provides another, without requiring connected marginal supports.

\begin{proposition}[Strict convexity]\label{prop:strictunique}
Under Assumption \ref{ass:base}, if $\psi$ is strictly convex on $\mathbb R$, then Assumption \ref{ass:unique} holds.
\end{proposition}
\begin{proof}
Let $h_i=f_i\oplus g_i$, $i=1,2$, be two canonical population optimal sums. If $h_1$ and $h_2$ differ on a set of positive $\mu\otimes\nu$ measure, strict convexity of $\psi$ makes
\begin{equation*}
 h\longmapsto \int\{h-\psi(h-c)\}\,\mathrm d(\mu\otimes\nu)
\end{equation*}
strictly concave along the segment joining $h_1$ and $h_2$, contradicting optimality. Thus $h_1=h_2$ $\mu\otimes\nu$-almost surely. Continuity extends the equality to $S_X\times S_Y$, which implies that $f_1-f_2$ and $g_1-g_2$ are opposite constants on the supports. The canonical extension equations then propagate the same gauge relation off the supports, and \eqref{eq:populationnormalization} fixes the gauge.
\end{proof}

A second curvature-based criterion can be stated directly in terms of the regularization strength. Restoring the parameter $\epsilon>0$, define the directional oscillations on the population supports by
\begin{align}\label{eq:directionalosc}
 \Delta_X(c)&:=\sup_{x,x'\in S_X}\sup_{y\in S_Y}|c(x,y)-c(x',y)|,\nonumber\\
 \Delta_Y(c)&:=\sup_{y,y'\in S_Y}\sup_{x\in S_X}|c(x,y)-c(x,y')|.
\end{align}
They are finite because $c$ is bounded.

\begin{proposition}[Large-$\epsilon$ uniqueness]\label{prop:largeeps}
Under Assumption \ref{ass:base}, the population potentials for $S_\epsilon$ are unique up to gauge whenever
\begin{equation}\label{eq:largeeps}
 \epsilon>\frac{2}{\delta}\min\{\Delta_X(c),\Delta_Y(c)\}.
\end{equation}
No connectedness assumption on the marginal supports is needed.
\end{proposition}
\begin{proof}
Scale to unit regularization with $\bar c=c/\epsilon$ and let $(\bar f,\bar g)$ be canonical potentials. Write
\begin{equation*}
 \xi(x,y)=\bar f(x)+\bar g(y)-\bar c(x,y).
\end{equation*}
By \eqref{eq:directionalg}, $\operatorname{osc}_{S_Y}(\bar g)\le\Delta_Y(c)/\epsilon$. Hence, for every fixed $x\in S_X$,
\begin{equation*}
 \operatorname{osc}_{y\in S_Y}\xi(x,y)\le\frac{2\Delta_Y(c)}{\epsilon}.
\end{equation*}
The first marginal equation $\int\psi'(\xi(x,y))\nu(\mathrm dy)=1=\psi'(t_0)$ forces $\sup_{y\in S_Y}\xi(x,y)\ge t_0$; otherwise monotonicity, together with strict increase of $\psi'$ on $[t_0-\delta,t_0]$, would make the integral strictly smaller than one. Combined with the oscillation estimate, this gives
\begin{equation}\label{eq:xilowerY}
 \inf_{y\in S_Y}\xi(x,y)\ge t_0-\frac{2\Delta_Y(c)}{\epsilon}.
\end{equation}
The symmetric argument gives the same bound with $\Delta_X(c)$. If \eqref{eq:largeeps} holds, one of these bounds is strictly larger than $t_0-\delta$, so the entire product $S_X\times S_Y$ lies in the strict-curvature region. Thus $\mathcal A^*=S_X\times S_Y$ for the scaled problem. Equality $u\oplus v=0$ on the full product immediately implies $u\equiv a$ and $v\equiv-a$, and Proposition \ref{prop:activeunique} completes the proof.
\end{proof}

\begin{remark}[Polynomial regularization]\label{rem:polynomialactive}
For
\begin{equation*}
 \phi_p(r)=\frac{r^p-1}{p},\qquad 1<p\le2,
\end{equation*}
let $q=p/(p-1)$. Then
\begin{equation*}
 \psi_p(s)=\frac{(s_+)^q}{q}+\frac1p,
 \qquad
 \psi_p'(s)=(s_+)^{q-1}.
\end{equation*}
One may take $t_0=1$ and $\delta=1$, so $r_\delta=0$ and $\mathcal A^*$ is exactly the positive-density set $\{\rho^*>0\}$. The case $p=2$ is quadratic regularization: $\psi_2(s)=\tfrac12(s_+)^2+\tfrac12$ is $C^1$ but not $C^2$ at zero. Thus Theorem \ref{thm:connectedunique} applies to this boundary case on arbitrary Polish spaces. In the finitely supported setting, Corollary \ref{cor:activegraph} reduces uniqueness to connectedness of the positive-support bipartite graph. Proposition \ref{prop:largeeps} becomes the explicit sufficient condition $\epsilon>2\min\{\Delta_X(c),\Delta_Y(c)\}$.
\end{remark}

We now turn to empirical stability. For this purpose, assume in addition that the cost is Lipschitz in each variable:
\begin{equation}\label{eq:lipcost}
 |c(x,y)-c(x',y)|\le Ld_x(x,x'),\qquad
 |c(x,y)-c(x,y')|\le Ld_y(y,y').
\end{equation}
Proposition \ref{prop:dualbasics} then gives uniform boundedness and a common Lipschitz constant for the canonical empirical potentials. Because empirical optimizers need not be unique, we first fix measurable canonical selections.

\begin{lemma}[Measurable selection]\label{lem:measselection}
Fix $n,m\ge1$. There exist normalized canonical dual optimizers that can be chosen measurably for each of the problems $S(\mu_n,\nu_m)$, $S(\mu_n,\nu)$, and $S(\mu,\nu_m)$. For the two-sample problem, in particular, the maps from the sample and evaluation point to $f_{n,m}(x)$ and $g_{n,m}(y)$ are jointly Borel measurable, and the selected pair satisfies
\begin{equation*}
 \mu_n f_{n,m}=0,
 \qquad
 \|f_{n,m}\|_\infty+\|g_{n,m}\|_\infty\le B,
\end{equation*}
with the deterministic constant $B$ from Proposition \ref{prop:dualbasics}.
\end{lemma}
\begin{proof}
Write the sample as $z=(x_1,\ldots,x_n,y_1,\ldots,y_m)$. Only finitely many equality patterns can occur among the $x_i$'s and among the $y_j$'s, and each pattern defines a Borel stratum of the sample space. On a fixed stratum, represent a dual pair by vectors $u=(u_1,\ldots,u_n)$ and $v=(v_1,\ldots,v_m)$ satisfying
\begin{equation*}
 u_i=u_{i'}\ \text{whenever }x_i=x_{i'},
 \qquad
 v_j=v_{j'}\ \text{whenever }y_j=y_{j'},
 \qquad
 \frac1n\sum_{i=1}^n u_i=0,
\end{equation*}
and $|u_i|,|v_j|\le B_1$, where $B_1:=\sqrt{n+m}\,B+1$. On a fixed stratum these constraints define a nonempty compact set $K\subset\mathbb R^{n+m}$ independent of the particular sample point within that stratum. The empirical dual objective is
\begin{equation*}
 \mathcal L_z(u,v)
 =\frac1n\sum_{i=1}^n u_i+\frac1m\sum_{j=1}^m v_j
 -\frac1{nm}\sum_{i=1}^n\sum_{j=1}^m
 \psi\big(u_i+v_j-c(x_i,y_j)\big).
\end{equation*}
The map is continuous in $(u,v)$ and Borel measurable in $z$, and Proposition \ref{prop:dualbasics} ensures that its maximum over $K$ equals $S(\mu_n,\nu_m)$. Among the maximizers choose one with minimal Euclidean norm. The measurable maximum theorem, applied successively to the objective and to the norm on the argmax correspondence, gives a measurable choice on each stratum; the finitely many selections can then be pasted over the full sample space. Since a normalized canonical optimizer with coordinatewise bound $B$ exists by Proposition \ref{prop:dualbasics}, the selected minimizer has Euclidean norm at most $\sqrt{n+m}\,B<B_1$. Hence none of the box constraints is active at the selected point.

It remains to verify the marginal equations on the empirical support. Let $\bar x_1,\ldots,\bar x_A$ and $\bar y_1,\ldots,\bar y_C$ be the distinct empirical atoms, with masses $p_a$ and $q_b$, and let $\bar u_a,\bar v_b$ be the corresponding selected potential values. In these coordinates the gauge is $\sum_a p_a\bar u_a=0$. Since the box constraints are inactive, differentiation with respect to $\bar v_b$ gives
\begin{equation*}
 \sum_{a=1}^A p_a\,\psi'\big(\bar u_a+\bar v_b-c(\bar x_a,\bar y_b)\big)=1,
 \qquad b=1,\ldots,C.
\end{equation*}
Introducing a Lagrange multiplier $\lambda$ for the gauge constraint, differentiation with respect to $\bar u_a$ gives
\begin{equation*}
 1-\sum_{b=1}^C q_b\,\psi'\big(\bar u_a+\bar v_b-c(\bar x_a,\bar y_b)\big)=\lambda,
 \qquad a=1,\ldots,A.
\end{equation*}
Multiplying by $p_a$, summing over $a$, and using the preceding column equations yields $\lambda=0$. Thus the row equations hold as well.

To extend the selected support values canonically, define for $x\in\mathcal X$
\begin{equation*}
 F_z(x,s):=\frac1m\sum_{j=1}^m
 \psi'\big(s+v_j-c(x,y_j)\big).
\end{equation*}
Dual regularity gives a unique finite solution of $F_z(x,s)=1$. It is jointly measurable in $(z,x)$ through the generalized-inverse representation
\begin{equation*}
 f_{n,m}(x)=\inf\{s\in\mathbb R:F_z(x,s)\ge1\},
\end{equation*}
and measurability follows because its strict sublevel sets can be written using rational test points and the joint measurability of $F_z(x,q)$. The analogous construction using the selected $u_i$'s defines $g_{n,m}(y)$ measurably. The roots are finite by monotonicity and the normalization $\psi'(t_0)=1$. The row and column equations established above show that these scalar roots reproduce the selected values on the empirical supports; hence the extension is canonical. Proposition \ref{prop:dualbasics} then gives the stated deterministic bound.

The one-sample problem $S(\mu_n,\nu)$ can be treated directly through the canonical equations. On a fixed equality-pattern stratum, let $u=(u_1,\ldots,u_n)$ satisfy the equality constraints for repeated sample points, the gauge $n^{-1}\sum_i u_i=0$, and $|u_i|\le B$. For fixed $(z,u)$ and $y\in\mathcal Y$, set
\begin{equation*}
 G_{z,u}(y,s)=\frac1n\sum_{i=1}^n
 \psi'\big(u_i+s-c(x_i,y)\big).
\end{equation*}
The equation $G_{z,u}(y,s)=1$ has a unique solution $g_{z,u}(y)$ in a deterministic bounded interval. The scalar-root argument used above shows that $(z,u,y)\mapsto g_{z,u}(y)$ is jointly Borel measurable and continuous in $u$ for fixed $(z,y)$. For each distinct empirical atom $\bar x_a$ with mass $p_a$, define the row residual
\begin{equation*}
 \mathcal R_a(z,u)
 :=\int \psi'\big(u(\bar x_a)+g_{z,u}(y)-c(\bar x_a,y)\big)\,\nu(\mathrm dy)-1.
\end{equation*}
The vector of residuals is continuous in $u$ and Borel measurable in $z$. Hence the set of gauge-normalized $u$ in the compact box $[-B,B]^n$ satisfying all equality constraints and $\mathcal R_a(z,u)=0$ for every distinct atom is a measurable compact-valued correspondence. It is nonempty by Proposition \ref{prop:dualbasics}. The Kuratowski--Ryll-Nardzewski measurable selection theorem gives a measurable choice $u(z)$. The associated $g_{z,u(z)}$ and the first scalar marginal equation then provide jointly measurable canonical extensions of both potentials. This yields a measurable canonical optimizer for $S(\mu_n,\nu)$. The construction for $S(\mu,\nu_m)$ is symmetric. Thus all empirical optimizers used below can be fixed as measurable canonical selections without assuming empirical uniqueness.
\end{proof}

We next establish stability for canonical empirical selections.

\begin{lemma}[Potential stability]\label{L2}
Assume \ref{ass:base}, \ref{ass:unique}, and \eqref{eq:lipcost}. Let $\mu_k\Rightarrow\mu$ and $\nu_k\Rightarrow\nu$ weakly, and let $(f_k,g_k)$ be any canonical optimal pair for $S(\mu_k,\nu_k)$ normalized by
\begin{equation*}
 \int f_k\,\mathrm d\mu_k=0.
\end{equation*}
Then
\begin{equation}\label{eq:localuniform}
 f_k\to f^*\quad\text{and}\quad g_k\to g^*
 \qquad\text{locally uniformly on }\mathcal X\text{ and }\mathcal Y.
\end{equation}
\end{lemma}
\begin{proof}
Proposition \ref{prop:dualbasics} gives deterministic uniform bounds and a common Lipschitz constant for $\{f_k\}$ and $\{g_k\}$. Fix countable dense subsets of $\mathcal X$ and $\mathcal Y$. From every subsequence, a diagonal argument extracts a further subsequence converging on both dense sets. Equi-Lipschitz continuity extends the pointwise limits to Lipschitz functions $(f_\infty,g_\infty)$ and upgrades the convergence to uniform convergence on compact sets.

Set $h_k=f_k\oplus g_k$ and $h_\infty=f_\infty\oplus g_\infty$. Since $h_k$ is uniformly bounded and converges locally uniformly, the same is true of $\ell_{h_k}=h_k-\psi(h_k-c)$. Weak convergence $\mu_k\otimes\nu_k\Rightarrow\mu\otimes\nu$, together with tightness and the uniform bound, yields
\begin{equation}\label{eq:movingintegral}
 L_{\mu_k,\nu_k}(h_k)\longrightarrow L_{\mu,\nu}(h_\infty).
\end{equation}
On a compact set carrying arbitrarily large mass the integrands converge uniformly, while the common bound controls the complement. The same argument gives
\begin{equation*}
 L_{\mu_k,\nu_k}(h^*)\longrightarrow L_{\mu,\nu}(h^*)=S(\mu,\nu).
\end{equation*}
Optimality of $h_k$ therefore implies
\begin{equation*}
 L_{\mu,\nu}(h_\infty)\ge L_{\mu,\nu}(h^*),
\end{equation*}
so $(f_\infty,g_\infty)$ is a population optimal pair. To verify that it is canonical, fix $x\in\mathcal X$. The first marginal equation, local uniform convergence, tightness of $\{\nu_k\}$, and the deterministic bound on the arguments of $\psi'$ give
\begin{equation*}
 1=\int\psi'(f_k(x)+g_k(y)-c(x,y))\,\nu_k(\mathrm dy)
 \longrightarrow
 \int\psi'(f_\infty(x)+g_\infty(y)-c(x,y))\,\nu(\mathrm dy).
\end{equation*}
The second marginal equation passes to the limit analogously for each fixed $y\in\mathcal Y$, so $(f_\infty,g_\infty)$ satisfies \eqref{eq:margx}--\eqref{eq:margy}. Applying the same moving-integral argument to $f_k$ gives $\int f_\infty\,\mathrm d\mu=0$. Assumption \ref{ass:unique} therefore implies $(f_\infty,g_\infty)=(f^*,g^*)$. Since every subsequence has a further subsequence with this limit, the full sequence converges as in \eqref{eq:localuniform}.
\end{proof}

For empirical marginals, local-uniform stability strengthens to the $L^2$ convergence needed in the remainder estimates.

\begin{corollary}[Empirical $L^2$ stability]\label{cor:l2stability}
Under the assumptions of Lemma \ref{L2}, almost surely,
\begin{align}
 \|h_n^{X}-h^*\|_{L^2(\mu_n\otimes\nu)}&\to0,\label{eq:l2x}\\
 \|h_m^{Y}-h^*\|_{L^2(\mu\otimes\nu_m)}&\to0,\label{eq:l2y}\\
 \|h_{n,m}-h^*\|_{L^2(\mu_n\otimes\nu_m)}&\to0\qquad(n,m\to\infty),\label{eq:l2xy}
\end{align}
where the superscripts indicate the one-sample problems. They also hold in mean square with respect to the sampling probability; equivalently, the expectations of the squared displayed norms converge to zero.
\end{corollary}
\begin{proof}
Varadarajan's theorem gives the required almost-sure weak convergence of the empirical measures, so Lemma \ref{L2} applies. For the double limit in \eqref{eq:l2xy}, work on the probability-one event on which the full sequences satisfy $\mu_n\Rightarrow\mu$ and $\nu_m\Rightarrow\nu$. If the two-index convergence failed on such a sample path, there would exist a (possibly path-dependent) sequence $(n_k,m_k)$ with $n_k,m_k\to\infty$ along which the norm stays bounded away from zero. Applying the deterministic Lemma \ref{L2} along that sequence gives a contradiction. Hence the stated two-index convergence holds almost surely.

To pass from local uniform convergence to the empirical $L^2$ norms, fix $\eta>0$ and choose compacts $K_X,K_Y$ with $\mu(K_X),\nu(K_Y)>1-\eta$. By the strong law applied to the indicators of the fixed compacts, $\mu_n(K_X)\to\mu(K_X)$ and $\nu_m(K_Y)\to\nu(K_Y)$ almost surely; hence, eventually, the corresponding empirical measures assign mass at least $1-2\eta$ to these compacts. On $K_X\times K_Y$ the convergence is uniform, while outside it the squared difference is bounded by a deterministic constant from Proposition \ref{prop:dualbasics}. Letting $\eta\downarrow0$ proves the almost-sure claims. Uniform boundedness then gives convergence of their expectations by dominated convergence.
\end{proof}

The leave-one-out argument below evaluates random potentials at sample-dependent points and integrates them against moving empirical measures. The next lemma handles both operations by a compact--tail argument.

\begin{lemma}\label{lem:movingrandom}
Let $q_k:\mathcal X\times\mathcal Y\to\mathbb R$ be jointly measurable random continuous functions such that, almost surely,
\begin{equation*}
 q_k\longrightarrow0\quad\text{locally uniformly on }\mathcal X\times\mathcal Y,
 \qquad
 \sup_k\|q_k\|_\infty\le C
\end{equation*}
for a deterministic $C<\infty$. Let $Z$ be an $\mathcal X$-valued random variable; no independence between $Z$ and $q_k$ is required.
\begin{enumerate}
 \item For every probability measure $Q$ on $\mathcal Y$,
 \begin{equation*}
  \int |q_k(Z,y)|\,Q(\mathrm dy)\longrightarrow0
  \qquad\text{almost surely and in }L^p
 \end{equation*}
 for every finite $p\ge1$.
 \item If $Q_k$ are random probability measures with $Q_k\Rightarrow Q$ almost surely, then
 \begin{equation*}
  \int |q_k(Z,y)|\,Q_k(\mathrm dy)\longrightarrow0
  \qquad\text{almost surely and in }L^p
 \end{equation*}
 for every finite $p\ge1$.
\end{enumerate}
The symmetric statements with the roles of $\mathcal X$ and $\mathcal Y$ interchanged also hold.
\end{lemma}
\begin{proof}
Fix an event of probability one on which all stated convergences hold, and write $z=Z(\omega)$. For the first assertion, given $\eta>0$, choose a compact $K\subset\mathcal Y$ with $Q(K)>1-\eta$. Since $\{z\}\times K$ is compact,
\begin{equation*}
 \sup_{y\in K}|q_k(z,y)|\longrightarrow0,
\end{equation*}
whereas the integral over $K^c$ is at most $C\eta$. Letting $\eta\downarrow0$ proves the almost-sure convergence.

For the second assertion, weak convergence $Q_k\Rightarrow Q$ implies tightness of $Q$ together with all sufficiently large $Q_k$. Hence, for every $\eta>0$, a compact $K$ can be chosen so that $Q_k(K^c)\le\eta$ eventually. The same compact--tail decomposition gives almost-sure convergence. In both cases the integrals are bounded by $C$, and dominated convergence yields convergence in every finite $L^p$.
\end{proof}

Define the bounded population functions
\begin{align}\label{eq:influence}
 \Psi^*(x,y)&:=\psi(h^*(x,y)-c(x,y)),\nonumber\\
 a(x)&:=f^*(x)-\int\Psi^*(x,y)\nu(\mathrm dy),\qquad
 b(y):=g^*(y)-\int\Psi^*(x,y)\mu(\mathrm dx).
\end{align}
The variances
\begin{equation}\label{eq:sigmabasic}
 \sigma_X^2:=\operatorname{Var}_\mu(a(X)),
 \qquad
 \sigma_Y^2:=\operatorname{Var}_\nu(b(Y))
\end{equation}
are invariant under the additive gauge and hence are well defined.

\begin{theorem}[Asymptotic linearization]\label{v1}
Assume \ref{ass:base}, \ref{ass:unique}, and \eqref{eq:lipcost}. Define
\begin{align*}
 T_n^X&:=\mu_n f^*+\nu g^*-(\mu_n\otimes\nu)\Psi^*,\\
 T_m^Y&:=\mu f^*+\nu_m g^*-(\mu\otimes\nu_m)\Psi^*,\\
 T_{n,m}&:=\mu_n f^*+\nu_m g^*-(\mu_n\otimes\nu_m)\Psi^*.
\end{align*}
By Fubini's theorem and, for $T_{n,m}$, independence of the two samples,
\begin{equation}\label{eq:Tunbiased}
 \mathrm E T_n^X=\mathrm E T_m^Y=\mathrm E T_{n,m}
 =\mu f^*+\nu g^*-(\mu\otimes\nu)\Psi^*=S.
\end{equation}
Define $R_n^X=S(\mu_n,\nu)-T_n^X$, $R_m^Y=S(\mu,\nu_m)-T_m^Y$, and $R_{n,m}=S(\mu_n,\nu_m)-T_{n,m}$. Then
\begin{align}
 n\operatorname{Var}(R_n^X)&\to0,\label{eq:remx}\\
 m\operatorname{Var}(R_m^Y)&\to0,\label{eq:remy}\\
 \frac{nm}{n+m}\operatorname{Var}(R_{n,m})&\to0\qquad(n,m\to\infty).\label{eq:remxy}
\end{align}
\end{theorem}
\begin{proof}
We prove \eqref{eq:remx} and \eqref{eq:remxy}; \eqref{eq:remy} follows symmetrically. Write $\ell_h=h-\psi(h-c)$. For the one-sample problem, let $h_n^X$ be an optimizer for $(\mu_n,\nu)$ and let $\mu_n'$ replace $X_1$ by an independent copy $X_1'$. Denote the corresponding remainder by $(R_n^X)'$. Using $h_n^X$ as a competitor in the primed problem gives
\begin{align*}
 R_n^X-(R_n^X)'
 &\le (\mu_n-\mu_n')\nu\big(\ell_{h_n^X}-\ell_{h^*}\big).
\end{align*}
If $(h_n^X)'$ denotes a canonical optimizer for the primed problem, exchanging the primed and unprimed samples gives the reverse estimate with $(h_n^X)'$. Consequently,
\begin{align}\label{eq:remdiffx}
 n\big|R_n^X-(R_n^X)'\big|
 &\le \left|\nu(\ell_{h_n^X}-\ell_{h^*})(X_1,\cdot)
 -\nu(\ell_{h_n^X}-\ell_{h^*})(X_1',\cdot)\right|\nonumber\\
 &\quad+\left|\nu(\ell_{(h_n^X)'}-\ell_{h^*})(X_1,\cdot)
 -\nu(\ell_{(h_n^X)'}-\ell_{h^*})(X_1',\cdot)\right|.
\end{align}
Lemma \ref{lem:measselection} makes all remainders and replacement differences measurable. Almost surely, $\mu_n\Rightarrow\mu$, and
\begin{equation*}
 \mu_n'=\mu_n+\frac1n(\delta_{X_1'}-\delta_{X_1})\Rightarrow\mu,
\end{equation*}
so Lemma \ref{L2} applies to both $h_n^X$ and $(h_n^X)'$, which therefore converge locally uniformly to $h^*$ almost surely. Since the map $h\mapsto \ell_h$ is uniformly continuous on the deterministic compact interval containing all $h-c$, the functions
\begin{equation*}
 \ell_{h_n^X}-\ell_{h^*}
 \quad\text{and}\quad
 \ell_{(h_n^X)'}-\ell_{h^*}
\end{equation*}
converge locally uniformly to zero and are uniformly bounded. Lemma \ref{lem:movingrandom}, applied with $Z=X_1$ and $Z=X_1'$ and $Q=\nu$, shows that all four terms on the right-hand side of \eqref{eq:remdiffx} converge to zero almost surely and in $L^2$. Hence
\begin{equation*}
 n^2\mathrm E\big|R_n^X-(R_n^X)'\big|^2\to0.
\end{equation*}
The standard Efron--Stein inequality now yields \eqref{eq:remx}.

For the two-sample problem, replace $X_1$ and $Y_1$ separately and write the corresponding remainders as $R_{n,m}^{X'}$ and $R_{n,m}^{Y'}$. To establish the two-index limits, fix an arbitrary deterministic sequence $(n_k,m_k)$ with $n_k,m_k\to\infty$. Then
\begin{equation*}
 \mu_{n_k},\ \mu_{n_k}^{X'}\Rightarrow\mu,
 \qquad
 \nu_{m_k},\ \nu_{m_k}^{Y'}\Rightarrow\nu
 \qquad\text{almost surely}.
\end{equation*}
Lemma \ref{L2} gives local-uniform convergence to $h^*$ for both the original and replaced canonical optimizers. Put
\begin{equation*}
 q_{n,m}:=\ell_{h_{n,m}}-\ell_{h^*},\qquad
 q_{n,m}^{X'}:=\ell_{h_{n,m}^{X'}}-\ell_{h^*},\qquad
 q_{n,m}^{Y'}:=\ell_{h_{n,m}^{Y'}}-\ell_{h^*}.
\end{equation*}
The same two-sided comparison used in \eqref{eq:remdiffx} gives
\begin{align*}
 n\big|R_{n,m}-R_{n,m}^{X'}\big|
 &\le \big|\nu_m q_{n,m}(X_1,\cdot)-\nu_m q_{n,m}(X_1',\cdot)\big|\\
 &\quad+\big|\nu_m q_{n,m}^{X'}(X_1,\cdot)-\nu_m q_{n,m}^{X'}(X_1',\cdot)\big|,
\end{align*}
and, symmetrically,
\begin{align*}
 m\big|R_{n,m}-R_{n,m}^{Y'}\big|
 &\le \big|\mu_n q_{n,m}(\cdot,Y_1)-\mu_n q_{n,m}(\cdot,Y_1')\big|\\
 &\quad+\big|\mu_n q_{n,m}^{Y'}(\cdot,Y_1)-\mu_n q_{n,m}^{Y'}(\cdot,Y_1')\big|.
\end{align*}
Along the sequence $(n_k,m_k)$, Lemma \ref{lem:movingrandom}, with $Q_k=\nu_{m_k}$ or $Q_k=\mu_{n_k}$, makes the right-hand sides converge to zero almost surely and in $L^2$. Thus, along every such sequence,
\begin{align*}
 n_k^2\mathrm E\big|R_{n_k,m_k}-R_{n_k,m_k}^{X'}\big|^2&\to0,\\
 m_k^2\mathrm E\big|R_{n_k,m_k}-R_{n_k,m_k}^{Y'}\big|^2&\to0.
\end{align*}
Thus the corresponding two-index limits hold as $n,m\to\infty$. Efron--Stein gives
\begin{equation*}
 \operatorname{Var}(R_{n,m})
 \le \frac12\left[
 n\mathrm E\big|R_{n,m}-R_{n,m}^{X'}\big|^2
 +m\mathrm E\big|R_{n,m}-R_{n,m}^{Y'}\big|^2\right].
\end{equation*}
Multiplying by $nm/(n+m)$ proves \eqref{eq:remxy}.
\end{proof}

Theorem \ref{v1} reduces the centered fluctuation to empirical averages of the bounded influence functions $a$ and $b$; after removing the degenerate two-sample term, the Gaussian limit follows from the classical central limit theorem.

\begin{theorem}[Expectation-centered CLTs]\label{CLT1}
Under Assumptions \ref{ass:base}, \ref{ass:unique}, and \eqref{eq:lipcost},
\begin{align}
 \sqrt n\big(S(\mu_n,\nu)-\mathrm ES(\mu_n,\nu)\big)
 &\xrightarrow{D}\mathcal N(0,\sigma_X^2),\label{eq:cltx}\\
 \sqrt m\big(S(\mu,\nu_m)-\mathrm ES(\mu,\nu_m)\big)
 &\xrightarrow{D}\mathcal N(0,\sigma_Y^2).\label{eq:clty}
\end{align}
If, in addition, $m/(n+m)\to\lambda\in(0,1)$, then
\begin{equation}\label{eq:cltxy}
 \sqrt{\frac{nm}{n+m}}
 \big(S(\mu_n,\nu_m)-\mathrm ES(\mu_n,\nu_m)\big)
 \xrightarrow{D}
 \mathcal N\big(0,\lambda\sigma_X^2+(1-\lambda)\sigma_Y^2\big).
\end{equation}
Here $\mathcal N(0,0)$ is understood as the point mass at zero. Thus the first and second one-sample limits are nondegenerate when $\sigma_X^2>0$ and $\sigma_Y^2>0$, respectively, while the two-sample limit is nondegenerate when $\lambda\sigma_X^2+(1-\lambda)\sigma_Y^2>0$.
\end{theorem}
\begin{proof}
For the first one-sample problem,
\begin{equation*}
 T_n^X=\nu g^*+\mu_n a,
\end{equation*}
so
\begin{equation*}
 \sqrt n\big(T_n^X-\mathrm ET_n^X\big)
 =\frac1{\sqrt n}\sum_{i=1}^n\big(a(X_i)-\mu a\big)
 \xrightarrow{D}\mathcal N(0,\sigma_X^2)
\end{equation*}
by the classical i.i.d. central limit theorem. Theorem \ref{v1} implies
\begin{equation*}
 \sqrt n\big[(R_n^X-\mathrm ER_n^X)\big]\to0
 \quad\text{in }L^2,
\end{equation*}
so \eqref{eq:cltx} follows by Slutsky's theorem. The proof of \eqref{eq:clty} is symmetric.

For the two-sample problem, an exact Hoeffding decomposition gives
\begin{equation}\label{eq:lineardecomp}
 T_{n,m}-S
 =(\mu_n-\mu)a+(\nu_m-\nu)b
 -(\mu_n-\mu)(\nu_m-\nu)\Psi^*.
\end{equation}
Let
\begin{equation*}
 \widetilde\Psi(x,y)=\Psi^*(x,y)-\nu\Psi^*(x,\cdot)-\mu\Psi^*(\cdot,y)+(\mu\otimes\nu)\Psi^*.
\end{equation*}
Then the last term in \eqref{eq:lineardecomp} equals $(\mu_n\otimes\nu_m)\widetilde\Psi$. The kernel is degenerate in each coordinate, and independence gives
\begin{equation*}
 \operatorname{Var}\big((\mu_n\otimes\nu_m)\widetilde\Psi\big)
 =\frac{\|\widetilde\Psi\|_{L^2(\mu\otimes\nu)}^2}{nm}.
\end{equation*}
Thus this term is $o_P(\sqrt{(n+m)/(nm)})$. The first two terms in \eqref{eq:lineardecomp} are independent, and
\begin{align*}
 \sqrt{\frac{nm}{n+m}}(\mu_n-\mu)a
 &=\sqrt{\frac{m}{n+m}}\frac1{\sqrt n}\sum_{i=1}^n(a(X_i)-\mu a),\\
 \sqrt{\frac{nm}{n+m}}(\nu_m-\nu)b
 &=\sqrt{\frac{n}{n+m}}\frac1{\sqrt m}\sum_{j=1}^m(b(Y_j)-\nu b).
\end{align*}
The joint classical CLT therefore gives the Gaussian limit in \eqref{eq:cltxy} for $T_{n,m}$. The remainder is negligible after centering by Theorem \ref{v1}, so Slutsky's theorem completes the proof.
\end{proof}

To replace expectation centering by the population value, it remains to show that the bias is little-$o$ on the central-limit scale.

\begin{lemma}[Negligible bias]\label{lem:fastbias}
Under the assumptions of Theorem \ref{CLT1},
\begin{align}
 \sqrt n\big(\mathrm ES(\mu_n,\nu)-S\big)&\to0,\label{eq:fastbiasx}\\
 \sqrt m\big(\mathrm ES(\mu,\nu_m)-S\big)&\to0.\label{eq:fastbiasy}
\end{align}
If $m/(n+m)\to\lambda\in(0,1)$, then
\begin{equation}\label{eq:fastbiasxy}
 \sqrt{\frac{nm}{n+m}}\big(\mathrm ES(\mu_n,\nu_m)-S\big)\to0.
\end{equation}
\end{lemma}
\begin{proof}
The interpolation bound \eqref{eq:biaspointwise} is now used without discarding the factor measuring the distance between empirical and population dual sums. For $(\mu_n,\nu)$, the population marginal equation gives $\nu(p^*-1)=0$, and Lemma \ref{L31} yields
\begin{equation*}
 0\le R_n^X
 \le \|h_n^X-h^*\|_{L^2(\mu_n\otimes\nu)}
 \|\mu_n(p^*-1)\|_{L^2(\nu)}.
\end{equation*}
Consequently,
\begin{equation*}
 \sqrt n\,\mathrm ER_n^X
 \le
 \left(\mathrm E\|h_n^X-h^*\|_{L^2(\mu_n\otimes\nu)}^2\right)^{1/2}
 \left(n\,\mathrm E\|\mu_n(p^*-1)\|_{L^2(\nu)}^2\right)^{1/2}.
\end{equation*}
Corollary \ref{cor:l2stability} makes the first factor tend to zero, while the conditional-variance calculation from Theorem \ref{KK1} bounds the second uniformly. Since $\mathrm ER_n^X=\mathrm ES(\mu_n,\nu)-S$, this gives \eqref{eq:fastbiasx}; \eqref{eq:fastbiasy} follows symmetrically.

For two samples, Lemma \ref{L31} gives
\begin{equation*}
 R_{n,m}
 \le \|h_{n,m}-h^*\|_{L^2(\mu_n\otimes\nu_m)}
 \left(
 \|\nu_m(p^*-1)\|_{L^2(\mu_n)}^2
 +\|\mu_n(p^*-1)\|_{L^2(\nu_m)}^2
 \right)^{1/2}.
\end{equation*}
Cauchy--Schwarz, Corollary \ref{cor:l2stability}, and the estimates $C/m$ and $C/n$ from Theorem \ref{KK1} show
\begin{equation*}
 \mathrm ER_{n,m}=o\!\left(\sqrt{\frac1n+\frac1m}\right).
\end{equation*}
Multiplication by $\sqrt{nm/(n+m)}$ proves \eqref{eq:fastbiasxy}.
\end{proof}

\begin{theorem}[Population-centered CLTs]\label{LLL}
Under the assumptions of Theorem \ref{CLT1},
\begin{align*}
 \sqrt n\big(S(\mu_n,\nu)-S(\mu,\nu)\big)&\xrightarrow{D}\mathcal N(0,\sigma_X^2),\\
 \sqrt m\big(S(\mu,\nu_m)-S(\mu,\nu)\big)&\xrightarrow{D}\mathcal N(0,\sigma_Y^2).
\end{align*}
If $m/(n+m)\to\lambda\in(0,1)$, then
\begin{equation*}
 \sqrt{\frac{nm}{n+m}}
 \big(S(\mu_n,\nu_m)-S(\mu,\nu)\big)
 \xrightarrow{D}
 \mathcal N\big(0,\lambda\sigma_X^2+(1-\lambda)\sigma_Y^2\big).
\end{equation*}
\end{theorem}
\begin{proof}
Combine the expectation-centered limits of Theorem \ref{CLT1} with the negligible-bias statements of Lemma \ref{lem:fastbias} and use Slutsky's theorem.
\end{proof}

\begin{corollary}[Connected-support CLT]\label{cor:connectedclt}
Assume \ref{ass:base} and \eqref{eq:lipcost}. If either $\operatorname{supp}(\mu)$ or $\operatorname{supp}(\nu)$ is connected, then the expectation-centered conclusions of Theorem \ref{CLT1} and the population-centered conclusions of Theorem \ref{LLL} hold without separately imposing Assumption \ref{ass:unique}. In particular, no compactness of the supports and no finite-dimensional structure of $\mathcal X$ or $\mathcal Y$ are required.
\end{corollary}
\begin{proof}
Theorem \ref{thm:connectedunique} yields Assumption \ref{ass:unique}; the claim then follows from Theorems \ref{CLT1} and \ref{LLL}.
\end{proof}

\begin{example}[Infinite-dimensional Hilbert space]
Let $H$ be a separable infinite-dimensional Hilbert space, take $\mathcal X=\mathcal Y=H$, and consider the bounded Lipschitz cost
\begin{equation*}
 c(x,y)=\min\{\|x-y\|_H,1\}.
\end{equation*}
Suppose that one marginal has connected support. For instance, let it be a centered Gaussian measure $\mathcal N(0,Q)$ with positive, injective, trace-class covariance operator $Q$. Then
\begin{equation*}
 \operatorname{supp}\mathcal N(0,Q)
 =\overline{\operatorname{Ran}(Q^{1/2})}=H,
\end{equation*}
so the support is connected. Theorem \ref{thm:connectedunique} gives population uniqueness, and Corollary \ref{cor:connectedclt} applies to every divergence covered by Assumption \ref{ass:base}, including quadratic regularization.
\end{example}

\begin{remark}\label{rem:quadratic}
Unlike the compact Euclidean $C^2$-conjugate setting of \citet{GonzalezSanzEcksteinNutz2025}, the scalar-value CLT above is formulated for Polish spaces with a dual-regular $C^1$ conjugate. Empirical potential uniqueness is not required: Lemma \ref{lem:measselection} provides measurable canonical selections, while population identifiability determines their subsequential limits.
\end{remark}

\section{Conclusion}\label{Se5}
Finite-sample and asymptotic results rely on different properties of the dual problem. Finite-sample estimation of the transport value is controlled by bounded dual geometry: dual interpolation and additive Hoeffding projection yield the bias bound, while bounded differences control the stochastic fluctuation. By contrast, the asymptotic distribution depends on the stability of the population dual representation. Once population identifiability is available, the empirical potentials converge locally uniformly, the value becomes asymptotically linear, and the remaining bias is negligible on the central-limit scale.

The identifiability results show that dual uniqueness is governed by the region on which the regularizer identifies the additive potential sum. The strong active set $\mathcal A^*$ makes this mechanism explicit. Connectivity of the active region propagates equality of competing potentials, which explains the connected-support criterion on general Polish spaces and the active-graph criterion for finite marginals. The large-$\epsilon$ condition provides a complementary regime in which the entire product support lies in the strict-curvature region. In particular, these arguments accommodate quadratic regularization despite the lack of $C^2$ smoothness of its conjugate.

Our results assume bounded costs and concern the scalar regularized transport value. Within this scope they allow noncompact and infinite-dimensional Polish state spaces and dual-regular $C^1$ conjugates. Extending the analysis to unbounded costs would require moment or tail assumptions in place of uniform dual bounds. Quantitative rates for potential stability would sharpen the remainder estimates, while coupling- and potential-level limit theory would extend the asymptotic results beyond the scalar value.

\begin{appendices}
\section{Proof of the projection lemma}\label{app}
\begin{proof}[Proof of Lemma \ref{L31}]
Let $P,Q$ be probability measures and let $k\in L^2(P\otimes Q)$. Write
\begin{equation*}
 k_0=(P\otimes Q)k,
 \qquad
 k_1(x)=Qk(x)-k_0,
 \qquad
 k_2(y)=Pk(y)-k_0,
\end{equation*}
and
\begin{equation*}
 k_3(x,y)=k(x,y)-k_0-k_1(x)-k_2(y).
\end{equation*}
Then $k_0,k_1,k_2,k_3$ are mutually orthogonal in $L^2(P\otimes Q)$, and $k_3$ is orthogonal to every additive function $h=f\oplus g$. Hence
\begin{equation*}
 \langle h,k\rangle
 =\langle h,k_0+k_1+k_2\rangle.
\end{equation*}
By Cauchy--Schwarz,
\begin{align*}
 |\langle h,k\rangle|^2
 &\le \|h\|_{L^2(P\otimes Q)}^2
 \|k_0+k_1+k_2\|_{L^2(P\otimes Q)}^2\\
 &=\|h\|_{L^2(P\otimes Q)}^2
 \big(k_0^2+\|k_1\|_{L^2(P)}^2+\|k_2\|_{L^2(Q)}^2\big).
\end{align*}
Moreover,
\begin{equation*}
 \|Qk\|_{L^2(P)}^2=k_0^2+\|k_1\|_{L^2(P)}^2,
 \qquad
 \|Pk\|_{L^2(Q)}^2=k_0^2+\|k_2\|_{L^2(Q)}^2.
\end{equation*}
Therefore
\begin{equation*}
 \|k_0+k_1+k_2\|_2^2
 \le \|Qk\|_{L^2(P)}^2+\|Pk\|_{L^2(Q)}^2,
\end{equation*}
which is precisely the projection estimate \eqref{eq:additiveprojection}.
\end{proof}
\end{appendices}

\backmatter

\bmhead{Acknowledgements}
This research is supported by National Key R\&D Program of China (2024YFA1012401), the Science and Technology Commission of Shanghai Municipality (23JC1400501), and Natural Science Foundation of China (12241103).

\bibliography{document}

\end{document}